\documentclass[intlimits,a4paper,twoside,10pt]{amsart}
\usepackage{amsfonts,amssymb,amsmath,amsthm,mathtools,dsfont}
\usepackage[english]{babel}
\usepackage{hyperref}
\hypersetup{pdfauthor={RAL,GS,PT},pdftitle={Symmetries of fPDEs},%
            colorlinks, linktocpage=true, pdfstartpage=3, pdfstartview=FitV,%
    breaklinks=true, pdfpagemode=UseNone, pageanchor=true, pdfpagemode=UseOutlines,%
    plainpages=false, bookmarksnumbered, bookmarksopen=true, bookmarksopenlevel=1,%
    hypertexnames=true, pdfhighlight=/O,%
    urlcolor=blue, linkcolor=blue, citecolor=red,
        }

\numberwithin{equation}{section}

\DeclareMathOperator{\fder}{\mathcal{D}}
\DeclareMathOperator{\dd}{d}
\newcommand{\RRfder}[3][]{\sideset{_{#2}}{_{#3}^{#1}}\fder}

\DeclareMathOperator{\DD}{D}
\DeclareMathOperator{\e}{e}
\DeclareMathOperator{\pr}{\mathrm{pr}}
\newcommand{\media}[1]{{\langle #1\rangle}}  
\newcommand{\calcin}[2]{\left.#1\right|_{#2}}

\newcommand{\R}{{\mathds{R}}}
\newcommand{\Z}{{\mathds{Z}}}
\newcommand{\N}{{\mathds{N}}}
\newcommand{\M}{{\mathcal{M}}}

\DeclareSymbolFont{bbold}{U}{bbold}{m}{n}
\DeclareSymbolFontAlphabet{\mathbbold}{bbold}

\theoremstyle{definition}
\newtheorem{defRL}{Definition}[section]
\newtheorem{defRLpartial}[defRL]{Definition}

\DeclareFontEncoding{FML}{}{}\DeclareFontSubstitution{FML}{futm}{m}{it}
\DeclareSymbolFont{fourier}{FML}{futm}{m}{it}
\DeclareMathSymbol{\partialup}{\mathord}{fourier}{130}
\DeclareMathOperator{\fpartial}{\partial}
\newcommand{\rrfder}[3][]{\sideset{_{#2}}{_{#3}^{#1}}\fpartial}

\theoremstyle{plain}
\newtheorem{invariancefpde}{Theorem}[section]
\newtheorem{Lie}[invariancefpde]{Theorem}
\newtheorem{genprothm}[invariancefpde]{Theorem}

\newtheorem{coroll1}[invariancefpde]{Corollary}
\newtheorem{coroll2}[invariancefpde]{Corollary}

\theoremstyle{remark}

\theoremstyle{definition}
\newtheorem{grupposimmetria}{Definition}[section]

\newtheorem{infjet}[grupposimmetria]{Definition}

\usepackage{lmodern}
\usepackage{newtxtext}\usepackage{newtxmath}
\usepackage[cal=boondoxo]{mathalfa}
\usepackage[T1]{fontenc}

\begin{document}

\begin{abstract}
We provide a general theoretical framework allowing us to extend the classical Lie theory for partial differential equations to the case of equations of fractional order. We propose a general prolongation formula for the study of Lie symmetries in the case of an arbitrary finite number of independent variables. We also prove the Lie theorem in the case of fractional differential equations, showing that the proper space for the analysis of these symmetries is the infinite dimensional jet space.
\end{abstract}

\author{Rosario Antonio Leo}
\address[R.A. Leo]{Dipartimento di Matematica e Fisica ``Ennio De Giorgi'', Universit\`a del Salento, Via per Arnesano, 73100 -- Lecce, Italy}
\email{leora@le.infn.it}
\author{Gabriele Sicuro}
\address[G. Sicuro]{Centro Brasileiro de Pesquisas F\'isicas, Rua Dr. Xavier Sigaud, 150, 22290-180, Rio de Janeiro, Brazil}
\email{sicuro@cbpf.br}
\author{Piergiulio Tempesta}
\address[P. Tempesta]{Departamento de F\'{\i}sica Te\'{o}rica II (M\'{e}todos Matem\'{a}ticos de la f\'isica), Facultad de F\'{\i}sicas, Universidad Complutense de Madrid, 28040 -- Madrid, Spain and Instituto de Ciencias Matem\'aticas (CSIC-UAM-UC3M-UCM), C/ Nicol\'as Cabrera, No 13--15, 28049 Madrid, Spain}
\email{p.tempesta@fis.ucm.es,piergiulio.tempesta@icmat.es}
\dedicatory{This paper (in revised form) is now published in\\ \textup{Fract.~Calc.~Appl.~Anal., Vol.~20, No.~1 (2017), pp.~212--231,
\textsc{doi}: {10.1515/fca-2017-0011}},\\
and is available online at \url{http://www.degruyter.com/view/j/fca}.}
\title[Lie theory for fPDEs]{A foundational approach to the Lie theory for fractional order partial differential equations}

\date{\today}

\maketitle

\tableofcontents

\section{Introduction}
The purpose of this paper is to provide a general analytic framework for the construction of exact solutions of fractional partial differential equations (\textsc{fpde}s) by means of symmetry techniques. The study of \textsc{fpde}s has an intriguing history, dating back to XIX century, with the works of Abel, Riemann, Liouville, and many other mathematicians \cite{Ross1977}. In the last decades, there has been a resurgence of interest in \textsc{fpde}s, since these equations have been recognized as fundamental not only in pure mathematics, but also in statistical mechanics \cite{Barkai2000,Mandelbrot1968,Metzler1999}, economics \cite{Scalas2000}, engineering \cite{Sabatier2007}, and many other fields, for their manifold applications. Consequently, the search for exact solutions of both stationary and evolution \textsc{fpde}s is of great relevance. Although many analytic approaches have been largely (but not exhaustively) investigated, the determination of algebraic and geometric methods for the study of \textsc{fpde}s is essentially an open problem.

Concerning group-theoretical methods, in Ref.~\cite{BL} the invariance of a fractional diffusion equation under scaling transformations was established. Subsequently, in Refs.~\cite{gazizov2009,gazizov2011,LST} an extension of the Lie approach for standard differential equations to the case of \textsc{fpde}s was proposed. Inspired by these preliminary results, group theoretical methods have been very recently applied, for example, to fractional versions of the Harry--Dym equation \cite{Huang2014} and of the Burgers equation \cite{Wang2016}. A more general analysis of symmetries of time fractional evolution equations can be found in Ref.~\cite{Huang2015}. In all these works, only the case of equations involving fractional derivatives with respect to a single independent variable (typically the time variable) has been analyzed, keeping the other derivatives of integer order. The interest in group theoretical methods is based on the fact that, in the case of standard partial differential equations (\textsc{pde}s), the theory of Lie symmetries represents one of the most effective tools for the determination of classes of analytic solutions \cite{bluman,ibragimov1985,olver,vinogradovsymmetries}. Since the crucial observation, due to G.~Birkhoff, of the potential relevance of Lie theory in the study of the equations of fluid dynamics, in the middle of the XX century an impressive body of literature has been produced on the applications of symmetries to both linear and nonlinear problems. Lie theory for \textsc{pde}s is still a flourishing research area. Nevertheless, to our knowledge, a complete and rigorous mathematical theory for the case of \textsc{fpde}s is essentially lacking.

The present paper contributes to fill this gap, proposing a coherent approach to the notion of symmetries for \textsc{fpde}s, which generalizes the known theoretical results presented, in particular in Refs.~\cite{gazizov2009,gazizov2011}. We show the existence of a natural connection between the theory of Lie invariance of \textsc{fpde}s and the classical theory of Lie symmetries for \textsc{pde}s, once extended to the case of equations involving an infinite number of derivatives \cite{vinogradov}. In particular, our extension of group--theoretical methods to the case of \textsc{fpde}s strongly relies on the geometrical notions of infinite prolongation and infinite dimensional jet space $\mathcal{J}^\infty$ \cite{krasil2011,saunders,vinogradov}. Indeed, a \textsc{fpde} can also be interpreted as a system of \textsc{pde}s involving an infinite number of (classical) derivatives.

Our main results are the following. First, we prove a generalized prolongation formula for the action of a vector field on fractional operators involving an arbitrary number of independent variables. In the process, we recover as special cases the results already known to the literature. Second, by means of this prolongation formula, we derive the Lie condition for a vector field to be an infinitesimal generator of a symmetry group. 
As in the classical approach, the vector fields defined by the previous condition possess a twofold role. Indeed, they generate diffeomorphisms that foliate the manifold of the independent and dependent variables into group orbits, transforming solutions into solutions. More importantly, they allow a given \textsc{fpde} to be reduced into another one, possessing fewer independent variables, or a fractional ordinary differential equation (\textsc{fode}) to be reduced into an equation of lower order. Solving the reduced equations enables us to construct exact solutions of the original problem.

The results obtained in this paper pave the way to the creation of a complete theory of Lie symmetries for \textsc{fpde}s, that would include generalizations to higher order symmetries, nonclassical symmetries, master symmetries, etc.


The paper is organized as follows. In Section \ref{sec:1}, some basic definitions concerning Riemann--Liouville fractional calculus are reviewed. Section \ref{sec:2} is devoted to the notion of symmetry for \textsc{fpde}s. We show the correspondence between the fractional theory and the classical one and we explicitly derive a generalized prolongation formula expressing the action of the infinitely prolonged symmetry field on fractional operators. Finally, an application of our approach to a $(N+1)$-dimensional nonlinear \textsc{fpde} is presented.

\section{Fundamental definitions}\label{sec:1}
In this Section, we briefly review some basic definitions of fractional calculus. For an extensive and detailed treatment of the subject, and the proofs of the results sketched here, we refer to the monographs in Refs.~\cite{kilbas,samko}. 

Given $p\in\R$, we introduce the symbol
\begin{equation}\label{op}
\media{p}\coloneqq\begin{cases}[p]&p\geq 0,\\ -1& p<0.\end{cases}
\end{equation}
In the previous expression $[p]\in\Z$ is the floor function of $p$, $[p]\leq p<[p]+1$. We can now state the following fundamental
\begin{defRL}[Riemann--Liouville fractional operator] Let $p\in\R$ and $f\colon [a,b]\subset\R\to\R$ be an analytic function in $(a,b)$. We shall call
\begin{equation}
\RRfder[p]{a}{t}f(t)\coloneqq\frac{1}{\Gamma(1+\media{p}-p)}\frac{\dd^{\media{p}+1}}{\dd t^{\media{p}+1}}\int_a^t(t-\tau)^{\media{p}-p}f(\tau)\dd\tau,\qquad t\in(a,b),\label{R}
\end{equation}
the Riemann--Liouville (RL) fractional operator of order $p$ and terminals $(a,t)$. Here $\Gamma(x)\coloneqq\int_0^\infty t^{x-1}\e^{-t}\dd t$ is the Gamma function.
\end{defRL}
The operator defined above is usually called the \textit{RL fractional derivative} when $p>0$ and the \textit{RL fractional integral} for $p<0$. Observe that the RL operator is linear and, for $k\in\Z$, 
\begin{equation}
\lim_{p\to k^\pm}\RRfder[p]{a}{t}f(t)=\begin{cases}\frac{\dd^k f(t)}{\dd t^k}&\text{if $k>0$,}\\f(t)&\text{if $k=0$,}\\
\frac{1}{\Gamma(-k)}\int_a^t(t-\tau)^{|k|-1}f(\tau)\dd\tau&\text{if $k<0$.}\end{cases}
\end{equation}
Similarly, we can introduce the partial RL operator \cite{kilbas} as follows.
\begin{defRLpartial}\label{DefRLpartial} Let $p_i\in\R$, $i=1,\dots, N$, and \[f(x_1,\dots,x_N)\colon [a_1,b_1]\times\cdots\times[a_N,b_N]\subset\R^N\to\R,\] be an analytic function in $(a_1,b_1)\times\dots\times(a_N,b_N)$. We define the mixed RL operator as
\begin{multline}\label{RLPO}
\rrfder[p_1,\dots,p_N]{a_1,\dots,a_N}{x_1,\dots,x_N}f(x_1,\dots,x_N)\coloneqq\\=\left(\prod_{i=1}^N\frac{\partial^{\media{p_i}+1}}{\partial x_i^{\media{p_i}+1}}\int_{a_i}^{x_i}\!\dd\xi_i\,\frac{(x_i-\xi_i)^{\media{p_i}-p_i}}{\Gamma(1+\media{p_i}-p_i)}\right)f(\xi_1,\dots,\xi_N).
\end{multline}
\end{defRLpartial}
Remarkably, a mixed RL operator acting on an analytic function can be formally written as
\begin{equation}\label{prop}
\rrfder[p_1,\dots,p_N]{a_1,\dots,a_N}{x_1,\dots,x_N}=\prod_{i=1}^N\left[\sum_{k_i=0}^\infty\binom{p_i}{k_i}\frac{(x_i-a_i)^{k_i-p_i}}{\Gamma(k_i-p_i+1)}\frac{\partial^{k_i}}{\partial x_i^{k_i}}\right].
\end{equation}
In the previous expression we have adopted the shorthand notation
\begin{equation}
    \binom{p}{k}\coloneqq\frac{\prod_{i=1}^k(p-i+1)}{\Gamma(k+1)}.
\end{equation}
Formula \eqref{prop} will be fundamental to link the theory of symmetries of fractional differential equations with the classical theory of symmetries of \textsc{pde}s. It has been proved, for example, in Ref.~\cite[Sec.~3.5]{oldham1974fractional} in the case of a fractional derivative operator with respect to one variable only acting on a space of analytic functions. The case of many variables can be treated similarly.

\section{Lie Theory for Fractional \textsc{Pde}s}\label{sec:2}
In this Section we shall develop a theoretical framework allowing us to generalize the classical Lie theory to the case of fractional partial differential equations.
\subsection{Local Lie groups} We shall focus on \textsc{fpde}s with one dependent variable $u\in U\subseteq\R$ and $N$ independent variables $(x_1,\dots,x_N)\in X\coloneqq[a_1,+\infty)\times\dots\times[a_N,+\infty)\subseteq\R^N$, i.e., on equations of the form
\begin{equation}
\label{fpde}\mathcal E(x,u,\rrfder[\mathbf p]{\mathbf a}{x} u,\dots)=0,
\end{equation}
where $\mathcal{E}$ is a function involving a finite number of RL operators \begin{equation}
\rrfder[\mathbf p]{\mathbf a}{x}u\coloneqq\rrfder[p_1,\dots,p_N]{a_1,\dots,a_N}{{x_1},\dots,{x_N}}u(x),\qquad p_i\in\R\,\ \forall i.\label{derivataf}
\end{equation} \textit{We will assume that all RL operators acting upon the same variable $x_i$ have the same lower extreme $a_i$.} By analogy with the case of \textsc{pde}s we propose the following notion of symmetry group of a \textsc{fpde}.
\begin{grupposimmetria}[Symmetries] A \textit{continuous group $G$ of (point) symmetries} for the \textsc{fpde} in Eq.~\eqref{fpde} is a local group of diffeomorphisms that map solutions $(x,u)\in M\subset X\times U\eqqcolon\M$ (with $M$ open subset of $\M$) into solutions $g\cdot(x,u)=(\tilde x,\tilde u)=(\Xi_g(x,u),\Phi_g(x,u))\in \M$, $g\in G$, for some smooth functions $\Xi_g\colon M\to X$, $\Phi_g\colon M\to U$.\end{grupposimmetria}
We denote by $\mathbf v\in\mathfrak g$ an element of the Lie algebra $\mathfrak g$ associated with the Lie group of transformations $G$. We shall consider vector fields of the form
\begin{equation}\label{vTM}
\mathbf v=\sum_{i=1}^N\xi^i(x,u)\frac{\partial}{\partial x_i}+\phi(x,u)\frac{\partial}{\partial u}.
\end{equation}
Here, assuming that the local action of the symmetry group $G$ is given by \[(x,u)\xrightarrow{g}(\tilde x,\tilde u),\quad g\equiv g(\epsilon)=\e^{\epsilon\mathbf v},\quad\epsilon\in\R,\]
the following formulae hold:
\begin{equation}\label{xifi}
\xi^i(x,u)=\left.\frac{\dd\Xi^i_{g(\epsilon)}(x,u)}{\dd \epsilon}\right|_{\epsilon=0},\qquad \phi(x,u)=\left.\frac{\dd\Phi_{g(\epsilon)}(x,u)}{\dd\epsilon}\right|_{\epsilon=0}.
\end{equation}
If a RL operator of $u$ with respect to $x_i$ appears in $\mathcal E$, we impose the additional constraint 
\begin{equation}
a_i\xrightarrow{g}a_i \Rightarrow \xi^i(x_1,\dots,x_{i-1},a_i,x_{i+1},\dots x_N)\equiv 0\Rightarrow \rrfder[\mathbf p]{\mathbf a}{x}u(x)\xrightarrow{g}\rrfder[\mathbf p]{\mathbf a}{\tilde x}\tilde u(\tilde x).
\end{equation}
To evaluate the functions in Eqs.~\eqref{xifi}, we need to \textit{prolong} our function $u$ in a suitable extended space, as we will show below. 

\subsection{Jet space and general prolongation formula}To treat properly the problem of symmetries of \textsc{fpde}s, we need to implement an ``infinite'' prolongation procedure of vector fields over a suitable infinite-dimensional jet space $\mathcal J^\infty$. The theory of infinite dimensional jet bundles \cite{krasil2011,saunders,vinogradov} is a generalization of the classical finite dimensional theory. In particular, infinitely prolonged jet spaces were thoroughly investigated by A.~M.~Vinogradov \cite{vinogradov} in his study of symmetries of infinitely prolonged equations (higher symmetries).

In the following, we will use the standard multi-index notation for derivatives of a function $f(x)\equiv f(x_1,\dots,x_N)$. Precisely, for a set of indices $\sigma=(i_1,\dots,i_K)$, $1\leq i_k\leq N$, $k=1,\dots,K$, we shall introduce the notation \[\frac{\partial^{|\sigma|}f(x)}{\partial x^\sigma}\coloneqq\left(\prod_{k=1}^K\frac{\partial}{\partial x_{i_k}}\right)f(x),\quad\text{with }|\sigma|\coloneqq K.\]

\subsubsection{Jet space} Let us start with the following standard definition \cite{krasil2011}.
\begin{infjet}[Infinite jet manifold] Let $(\mathcal M,\pi,\mathcal X)$ be a fiber bundle with $\pi\colon\mathcal M\to \mathcal X$ projection from the total space $\mathcal M$ to the base space $\mathcal X$. We assume $\dim\mathcal M=\dim\mathcal X+1=N+1$. Then, denoting by $\Gamma_\eta(\pi)$ the set of all local sections of $\pi$ containing $\eta\in\mathcal X$, we say that two smooth sections $f,\tilde f\in\Gamma_\eta(\pi)$, are \textit{$k$-equivalent} in $\eta$ if
\begin{equation}
\left.\frac{\partial^{|\sigma|}f(x)}{\partial x^\sigma}\right|_{x=\eta}=\left.\frac{\partial^{|\sigma|}\tilde f(x)}{\partial x^\sigma}\right|_{x=\eta},\quad\text{for any $0\leq|\sigma|<k$,}
\end{equation}
supposing that a system local coordinates on $\mathcal X$ is given. We call the equivalence class of $f$ the \textit{$k$-jet} of $f$ in $\eta$, and we shall denote it by $[f]^k_\eta$. The $k$-jet manifold $\mathcal J^k(\pi)$ is by definition the set of these equivalence classes
\begin{equation}
\mathcal J^k(\pi)\coloneqq \left\{[f]^k_\eta \colon\eta\in\mathcal X,\ f\in\Gamma_\eta(\pi)\right\}.
\end{equation}
Moreover, the smooth fibre bundles
\begin{subequations}\begin{align}
\pi_k\colon& \mathcal J^k(\pi)\to \mathcal X,\quad [f]^k_\eta\mapsto\eta,\\
\pi_{k,l}\colon& \mathcal J^k(\pi)\to \mathcal J^l(\pi),\quad [f]^k_\eta\mapsto [f]^l_\eta,\quad k\geq l,
\end{align}\end{subequations}
are naturally introduced. The infinitely prolonged jet space $\mathcal J^\infty(\pi)$ is defined as the inverse limit of the sequence
\begin{equation}
\dots \to \mathcal J^{k+1}(\pi)\xrightarrow{\pi_{k+1,k}} \mathcal J^k(\pi)\to\dots\to \mathcal J^{1}(\pi)\xrightarrow{\pi_1}\mathcal X.
\end{equation}
\end{infjet}
In a similar manner, given a submanifold $\gamma\subseteq\M$, we denote by $\pr^k\gamma|_\eta$ the equivalence class of all submanifold having a $k$-th order contact with $\gamma$ in $\eta$, and $\pr^k\gamma\coloneqq\bigcup_\eta \pr^k\gamma|_\eta$. It is possible to define $\pr^\infty\gamma$ as a limit manifold obtained through a chain of projections using the maps $\pi_{k,l}$ introduced above. Finally, the prolongation of the action of a symmetry group $G$ acting on $\M$ is expressed similarly as $\pr^\infty g[\pr^\infty\gamma|_\eta]\coloneqq \pr^\infty [g\cdot \gamma|_\eta]$ for a submanifold $\gamma\subseteq\M$ and for $g\in G$. 

Given a fractional differential equation of the form of Eq.~\eqref{fpde}, the space introduced above represents a {suitable space for the prolongation procedure} we require, with the identifications $\mathcal M\equiv X\times U$ and $\mathcal X\equiv X$. Indeed, using Eq.~\eqref{prop}, we can rewrite Eq.~\eqref{fpde} as a differential equation involving an \textit{infinite} number of ordinary partial derivatives. For example, let us consider the fractional differential equation
\begin{equation}
    \RRfder[p]{a}{x}u(x)=f(x),\quad p\in\R.
\end{equation}
Then, we can write
\begin{equation}\label{sistema}
  \RRfder[p]{a}{x}u(x)=f(x)\Leftrightarrow 
  \begin{cases}
  z_0=f(x)\\
  \omega_0=u\\
  \omega_{k}=\frac{\dd \omega_{k-1}}{\dd x},&k\in\N,\\
  z_{k}=z_{k-1}-\binom{p}{k-1}\frac{(x-a)^{k-p-1}}{\Gamma(k-p)}\omega_{k-1},& k\in\N.
  \end{cases}
\end{equation}
The problem of local symmetries of a \textsc{fpde} can be solved therefore in the geometrical framework of the infinite dimensional jet space $\mathcal J^\infty$ of integer order derivatives introduced above; consequently, the theory of symmetry reduction is automatically borrowed from the classical theory. Eq.~\eqref{fpde} defines an infinite dimensional submanifold 
\begin{equation}
\Delta\coloneqq\left\{[f]^\infty_\eta\in \mathcal J^\infty(\pi)\colon\mathcal E(\eta,f,\rrfder[\mathbf p]{\mathbf a}{\eta} f,\dots)=0\right\}\subset\mathcal J^\infty(\pi).
\end{equation}
A solution of the given \textsc{fpde} is a smooth submanifold $\gamma\subseteq\M$, such that $\pr^\infty\gamma\subseteq\Delta$. An infinitesimal symmetry of the equation $\mathcal E$ is generated by a \textit{vector field} tangent to $\Delta$ which preserves the contact structure of $\Delta$. Observe that in the infinite dimensional case we cannot associate, in general, a \textit{flow} to $\mathbf v$ \cite{vinogradov}: if this happens, $\mathbf v$ is said to be a \textit{Lie field} (note that we deal with one dependent variable only). In the present paper, we restrict ourselves to Lie fields whose structure in local coordinates is
\begin{equation}
\mathbf v=\sum_{i=1}^N\xi^i(x,u)\frac{\partial}{\partial x_i}+\phi(x,u)\frac{\partial}{\partial u}+\sum_{|\sigma|\geq 1}\phi^\sigma(x,u,\dots,\partial^{\kappa}_xu,\dots)\frac{\partial^{|\sigma|}}{\partial u^\sigma},\quad u^\sigma\coloneqq\frac{\partial^{|\sigma|}u}{\partial x^\sigma},
\end{equation}
where $\xi^i$, $\phi$ and $\phi^\sigma$ are analytic functions of their variables. Introducing the \textit{total derivative operator}
\begin{equation}
\DD_i\coloneqq\frac{\partial}{\partial x_i}+\sum_{|\sigma|\geq 0}\frac{\partial u^\sigma}{\partial x_i}\frac{\partial}{\partial u^\sigma},
\end{equation}
it can be shown \cite{vinogradov} that \textit{every vector field associated to an infinitesimal symmetry transformation} can be written as
\begin{equation}\label{prv}
\pr^\infty\mathbf v=\sum_{i=1}^N\xi^i(x,u)\frac{\partial}{\partial x_i}+\phi(x,u)\frac{\partial}{\partial u}+\sum_{|\sigma|\geq 1}\phi^\sigma\frac{\partial}{\partial u^\sigma}
\end{equation}
where, for $\sigma\equiv(i_1,\dots,i_K)$, $\DD^\sigma\coloneqq \DD_{i_1}\circ\dots \circ\DD_{i_K}$, and
\begin{equation}\label{clprofor}
\phi^\sigma\coloneqq \DD^\sigma\left(\phi(x,u)-\sum_{i=1}^N\xi^i(x,u)\frac{\partial u}{\partial x_i}\right)+\sum_{i=1}^N\xi^i(x,u)\frac{\partial}{\partial x_i} \DD^\sigma u.
\end{equation}
The latter quantity depends on $N+1$ functions $\xi^1,\dots,\xi^N,\phi$ to be determined. We say that the vector field in Eq.~\eqref{prv} expresses the \textit{prolonged action} of the vector field
\begin{equation}
\mathbf v=\sum_{i=1}^N\xi^i(x,u)\frac{\partial}{\partial x_i}+\phi(x,u)\frac{\partial}{\partial u}\in\mathsf T M.
\end{equation}

Summarizing, in our perspective a \textsc{fpde} defines an infinite dimensional manifold in the space $\mathcal J^\infty$ endowed with an infinite order contact structure. The problem of finding a symmetry for a \textsc{fpde} (i.e., the corresponding functions $\xi^1,\dots,\xi^N,\phi$) therefore can be formulated on the same geometric structures introduced in the classical theory of higher symmetries \cite{vinogradov}. The problem can also be solved in the space $\mathcal J^\infty$ by imposing a suitable version of the standard Lie theorem that will be discussed below.

\subsubsection{General prolongation formula}
We look now for an explicit expression for the quantity
\begin{equation}
    \phi^\mathbf{p}\coloneqq \pr^{\infty}\mathbf v\left[\rrfder[\mathbf p]{\mathbf a}{x}u\right],
\end{equation}
as a function of the coefficients of $\textbf{v}$. To stress the formal analogy with Eq.~\eqref{clprofor}, we call the expression for $\phi^\mathbf{p}$ the \textit{general prolongation formula}. 

\begin{genprothm}[General prolongation formula]\label{GenProThm} Let $\mathbf v=\sum_{i=1}^N\xi^i(x,u)\partial_{x_i}+\phi(x,u)\partial_u$ be a vector field defined on $\mathsf T M$, $M\subseteq\M=X\times U$, where $X$ is an $N$--dimensional space, whose generic element has the form $x=(x_1,\dots,x_N)\in X$, and $U$ is a one dimensional space. The general prolongation formula is given by
\begin{equation}\label{genprolongation}
\phi^{\mathbf p}\coloneqq\pr^{\infty}\mathbf v\left[\rrfder[\mathbf p]{\mathbf a}{x}u\right]=\RRfder[\mathbf p]{\mathbf a}{x}\left(\phi-\sum_{i=1}^N\xi^i\partial_{x_i}u\right)+\sum_{i=1}^N\xi^i\partial_{x_i}\rrfder[\mathbf p]{\mathbf a}{x}u.\end{equation}
\end{genprothm}
\noindent In the previous equation we have introduced the fractional total derivative
\begin{equation}\label{RLPOT}
\RRfder[\mathbf p]{\mathbf a}{x}\coloneqq\prod_{i=1}^N\left[\sum_{k_i=0}^\infty\binom{p_i}{k_i}\frac{(x_i-a_i)^{k_i-p_i}}{\Gamma(k_i-p_i+1)}\DD_i^{k_i}\right].
\end{equation}
Note that the classical prolongation formula \cite{bluman,ibragimov1985,olver} is immediately recovered in the limit case of differentiation of integer order.
\begin{proof}
Using Eq.~\eqref{prop}, we shall apply the vector field in Eq.~\eqref{prv} to
\begin{equation}
\rrfder[\mathbf p]{\mathbf a}{x}u\coloneqq u\prod_{i=1}^N\frac{(x_i-a_i)^{-p_i}}{\Gamma(1-p_i)}+\sum_{|\sigma|\geq 1}C_\sigma(x)u^\sigma.
\end{equation}
In the latter expression, $\sigma=(\sigma_1,\dots,\sigma_K)$ is a multi-index and
\begin{equation}
C_\sigma(x)\coloneqq\prod_{i=1}^N\frac{(x_i-a_i)^{k_i(\sigma)-p_i}}{\Gamma\left(\sum_{j=1}^K\delta_{\sigma_j,i}-p_i+1\right)}.
\end{equation}
We obtain directly
\begin{multline}
\phi^{\mathbf p}\coloneqq\pr^{\infty}\mathbf v\left[\rrfder[\mathbf p]{\mathbf a}{x}u\right]=\\=\phi\prod_{i=1}^N\frac{(x_i-a_i)^{-p_i}}{\Gamma(1-p_i)}+\sum_{|\sigma|\geq 1}C_\sigma(x)\phi^\sigma
+u\sum_{i=1}^N\xi^i\frac{(x_i-a_i)^{-p_i-1}}{\Gamma(-p_i)}\prod_{j\neq i}\frac{(x_j-a_j)^{-p_j}}{\Gamma(1-p_j)}\\
+\sum_{i=1}^N\xi^i\sum_{|\sigma|\geq 1}u^\sigma\frac{(x_i-a_i)^{k_i(\sigma)-p_i-1}}{\Gamma(k_j(\sigma)-p_j)}\prod_{j\neq i}^N\frac{(x_j-a_j)^{k_j(\sigma)-p_j}}{\Gamma(k_j(\sigma)-p_j+1)}\\
=\RRfder[\mathbf p]{\mathbf a}{x}\phi-\sum_{i=1}^N\RRfder[\mathbf p]{\mathbf a}{x}\left(\xi^i\partial_{x_i}u\right)+\sum_{i=1}^N\xi^i\partial_{x_i}\rrfder[\mathbf p]{\mathbf a}{x}u,
\end{multline}
which coincides with the general prolongation formula in Eq.~\eqref{genprolongation}.
\end{proof}

\begin{coroll1}[Prolongation formula for $N=1$] \label{Coroll1} Under the hypotheses of Theorem \ref{GenProThm}, let us consider $N=1$ and $\mathbf v=\xi(x,u)\partial_x+\phi(x,u)\partial_u\in \mathsf T M$, $M\subset\M\coloneqq X\times U$. The general prolongation formula is given by
\begin{equation}\label{ord1}\phi^{p}
=\rrfder[p]{a}{x}\phi-\sum_{n=0}^\infty\binom{p}{n+1}\frac{\dd^{n+1}\xi}{\dd x^{n+1}}\rrfder[p-n]{a}{x}u.
\end{equation}
\end{coroll1}

For $a=0$ and $p>0$, it is convenient to make explicit the dependence on $u$ of the coefficient $\phi\equiv\phi(x,u(x))$ appearing in the vector field of Eq.~\eqref{ord1}. To this purpose, we shall use the following formula, due to Osler \cite{osler},
\begin{multline}
\RRfder[p]{0}{t}f(t,g(t))=\\=\sum\limits_{n=0}^\infty\sum\limits_{m=0}^n\sum\limits_{k=0}^m\sum\limits_{r=0}^k\binom{p}{n}\binom{n}{m}\binom{k}{r}\frac{t^{n-p}(-g)^r}{k!\Gamma(n+1-p)}\frac{\dd^mg^{k-r}}{\dd t^m}\frac{\partial^{n-m+k}f(t,g)}{\partial t^{n-m}\partial g^k}.\label{oslerformula}
\end{multline}
Isolating the terms that are linear in $u$ an its derivatives, we have
\begin{equation}\RRfder[p]{0}{x}\phi(x,u(x))=
\rrfder[p]{0}{x}\phi-u\rrfder[p]{0}{x}\left(\partial_u\phi\right)+\partial_u\phi\rrfder[p]{0}{x}u+\sum\limits_{k=1}^\infty\binom{p}{k}\frac{\partial^{k+1}\phi}{\partial x^k\partial u} \rrfder[p-k]{0}{x}u+\mu^p(x,u).\label{23}\end{equation}
In Eq.~\eqref{23} we have denoted by $\mu^p(x,u)$ the expression
\begin{equation}\mu^p(x,u)\coloneqq\sum\limits_{n=2}^\infty\sum\limits_{m=2}^n\sum\limits_{k=2}^m \sum\limits_{r=0}^{k-1}\binom{p}{n}\binom{n}{m}\binom{k}{r}\frac{x^{n-p}(-u)^r}{k!\Gamma(n+1-p)}\frac{\dd^m}{\dd x^m}\left(u^{k-r}\right)\frac{\partial^{n-m+k}\phi(x,u)}{\partial x^{n-m}\partial u^k}.\label{mu}\end{equation}
We deduce that, for $p>0$
\begin{multline}
\phi^p=\rrfder[p]{0}{x}\phi+\rrfder[p]{0}{x}u\left(\partial_u\phi-p\frac{\dd\xi}{\dd x}\right)
-u\rrfder[p]{0}{x}\partial_u\phi\\+\mu^p(x,u)+\sum\limits_{n=1}^\infty\left[\binom{p}{n}\partial^n_x\partial_u\phi
-\binom{p}{n+1}\frac{\dd^{n+1}\xi}{\dd x^{n+1}}\right]\rrfder[p-n]{0}{x}u.\label{prolpart}
\end{multline}
For further details on this particular case, we refer to Ref.~\cite{gazizov2011}.

\begin{coroll2} Assume that the hypotheses of Theorem \ref{GenProThm} hold with $N=2$, $x_1\equiv x$, $x_2\equiv t$ and $a_2=0$.  For $p>0$, we have
\begin{multline}\label{phi0pbis}\phi^p_t\coloneqq\phi^{0,p}=\rrfder[p]{0}{t}\phi+\rrfder[p]{0}{t}u\left(\partial_u\phi-p\frac{\dd\tau}{\dd t}\right)-u\rrfder[p]{0}{t}\partial_u\phi\\
+\sum\limits_{n=1}^\infty\left\{\left[\binom{p}{n}\partial^n_t\partial_u\phi-\binom{p}{n+1}\frac{\dd^{n+1}\tau}{\dd t^{n+1}}\right]\rrfder[p-n]{0}{t}u-\binom{p}{n}\frac{\dd^n\xi}{\dd t^n}\partial_x\rrfder[p-n]{0}{t}u\right\}\\
+\sum\limits_{n=2}^\infty\sum\limits_{m=2}^n\sum\limits_{k=2}^m\sum\limits_{r=0}^{k-1}\binom{p}{n}\binom{n}{m}\binom{k}{r}\frac{t^{n-p}(-u)^r}{k!\Gamma(n+1-p)}\frac{\dd^m}{\dd t^m}\left(u^{k-r}\right)\frac{\partial^{n-m+k}\phi}{\partial t^{n-m}\partial u^k}.
\end{multline}
A similar expression is obtained for $\phi^{p,0}\eqqcolon\phi^p_x$.
\end{coroll2}
The prolongation formula in Eq.~\eqref{phi0pbis} can be found in Ref.~\cite{LST}.

\subsection{The Lie theorem, its application and symmetry reduction}
We can now state the following fundamental result.
\begin{Lie}[Lie theorem for \textsc{fpde}s] Let $\mathcal E(x,u,\rrfder[p]{a}{x} u,\dots)=0$ be a \textsc{fpde}, defined on the open subset $M\subset\mathcal M$. If
\begin{equation}\label{lieformula}\calcin{\pr^{\infty}\mathbf v\left[\mathcal E\left(x,u,\rrfder[\mathbf p]{\mathbf a}{x} u,\dots\right)\right]}{\mathcal E=0}=0,
\end{equation}
then $\mathbf v=\sum_{i=1}^N\xi^i(x,u)\frac{\partial}{\partial x_i}+\phi(x,u)\frac{\partial}{\partial u}$ generates an infinitesimal symmetry of the equation $\mathcal E=0$.
\end{Lie}
Taking into account all the previous results, this theorem is a direct consequence of the Lie theorem stated in Ref.~\cite{vinogradov}, applied to the case of a \textsc{fpde} once it is written as a system of equations of infinite order, as done for example in Eq.~\eqref{sistema}. To apply the previous theorem to a \textsc{fpde}, we must (in principle) expand every RL operator using Eq.~\eqref{prop} and then apply the Lie condition using Theorem \ref{GenProThm} in Eq.~\eqref{lieformula} to determine the coefficients $\xi^i,\phi$ of all infinitesimal symmetry vector fields: the Lie equation provides an infinite number of differential equations for the $N+1$ components of the field $\mathbf v$. \textit{The symmetries of the general \textsc{fpde} are therefore obtained by solving the classical Lie condition in the standard infinite dimensional jet space $\mathcal J^\infty$}.

The determining equation coming from the Lie theorem can be quite involved, due to the infinite number of independent components in $\mathcal J^\infty$. Moreover, according to the previous analysis, different fractional operators cannot be considered, strictly speaking, to be independent in the space $\mathcal J^\infty$ (this fact can be easily seen by considering that the cardinality of the set of fractional operators is larger then the cardinality of the set of integer order derivatives). 
However, if we regard $\rrfder[\mathbf p]{\mathbf a}{x}u$ as an ``effective coordinate'' in $\mathcal J^\infty$ we can still obtain symmetries for the initial differential problem; in particular, the expressions for the action of the symmetry field on RL operators can be formulated in a more compact way. This fact, in turn, allows us to write down in a simpler form the Lie equation \eqref{lieformula} and the general prolongation formula in Eq.~\eqref{genprolongation}. Proceeding in this manner, we will find in general a \textit{subgroup} of the complete group of Lie symmetries of the equation. It may happen, however, that the specific form of the equation justifies this independence assumption. For example, let us suppose that we search for the solution of the following functional equation in $\mathcal J^\infty$:
\begin{equation}
\psi_1\left(x,u,\frac{\dd u}{\dd x},\dots,\frac{\dd^lu}{\dd x^l}\right)\RRfder[p]{a}{x}u(x)+\psi_2\left(x,u,\frac{\dd u}{\dd x},\dots,\frac{\dd^mu}{\dd x^m}\right)\RRfder[q]{a}{x}u(x)=0.
\end{equation}
with $p,q\in\R\setminus\N_0$ and $l,m\in\N$. The equation above can be written as
\begin{equation}
\sum_{k=0}^\infty\left[\binom{p}{k}\frac{(x-a)^{k-p}}{\Gamma(k-p+1)}\psi_1+\binom{q}{k}\frac{(x-a)^{k-q}}{\Gamma(k-q+1)}\psi_2\right]\frac{\partial^ku}{\partial x^k}=0.
\end{equation}
For $k\geq\max\{l,m\}$ we have
\begin{equation}
\psi_1\frac{\Gamma(p+1)(x-a)^{-p}}{\Gamma(p-k+1)\Gamma(k-p+1)}+\psi_2\frac{\Gamma(q+1)(x-a)^{-q}}{\Gamma(q-k+1)\Gamma(k-q+1)}=0.
\end{equation}
The fact that the previous relation must hold for all values of $k\geq\max\{l,m\}$, implies $\psi_1=\psi_2=0$.

\subsubsection{Reduction by symmetry} Given a local group $G$ of transformations of a \textsc{pde}, with corresponding Lie algebra $\mathfrak g$, the theory of symmetry reduction is formulated in terms of the quotient manifold  $\M/G$ of the orbits of $G$ on $\M$, and of the \textit{projection map} $\Pi\colon \M\to\M/G$ \cite{bluman,olver,vinogradov}. This map associates with every point $(x,u)\in\M$ its orbit under the action of $G$. In particular, the existence of the quotient manifold $\M/G$ and of the projection map $\Pi$ is assured if $G$ acts regularly, with $s$ dimensional orbits on the $(N+1)$ dimensional manifold $\M$. The linear map $\dd\Pi\colon \calcin{T\M}{(x,u)}\to\calcin{T(\M/G)}{\Pi(x,u)}$ is onto and its kernel is given by $\calcin{\mathfrak g}{(x,u)}\coloneqq\{\calcin{\mathbf v}{(x,u)}\colon\mathbf v\in\mathfrak g\}$. For a proof of the previous statements, see for example Ref.~\cite{olver}. In the context of \textsc{fpde}s, the same philosophy applies, because the symmetry algebra acts on the smooth manifold $\M$. For the sake of completeness, we sketch here the main results for our simple case. Additional details, and more refined derivations, can be found in the monograph by Olver \cite{olver}. 

On the basis of the existence of $\Pi$, we assume, as usual, that $\M/G \sim \R^{N+1-s}$ at least locally; therefore we can (locally) identify $N+1-s$ new functionally independent invariants $\{\eta^i(x,u)\}_i$ on $\M$ that parametrize $\M/G$. Among these new variables, we can choose a new dependent variable, let us call it $\upsilon(x,u)$, and $N-s$ new independent variables, say $\{\zeta^i(x,u)\}_{i}$. Moreover, the group $G$ identifies also an \textit{invariant infinite jet space} $\mathcal I^\infty(G)$ under the action of $G$ as the inverse limit of the chain
\begin{multline}
\dots \to \mathcal I^{k+1}(G)\xrightarrow{\pi_{k+1,k}} \mathcal I^k(G)\to\dots\to \mathcal I^{1}(G)\xrightarrow{\pi_1}\mathcal X,\\ \mathcal I^{k}(G)\coloneqq\left\{z\in \mathcal J^k\colon \exists\gamma\subseteq\M,\ \dim\gamma=N,\ \text{$G$-invariant such that $z=\left.\pr^{(k)}\gamma\right|_z$}\right\}.
\end{multline}
A (locally) $G$-invariant solution $\gamma\in\M$ of $\Delta$ is therefore such that $\pr^\infty\gamma\in\Delta\cap \mathcal I^\infty(G)$. Let us define now the projection $\Pi^\infty$ such that, for $z\in\M$ and $\gamma\subseteq\M$ $G$-invariant, $\Pi^\infty[\pr^\infty\gamma|_z]\coloneqq \pr^\infty(\gamma/G)|_{\Pi(z)}$. If $\gamma$ is a solution of the considered \textsc{fpde}, then $\pr^\infty(\gamma/G)=\Pi^\infty[\pr^\infty\gamma]\in \Pi^\infty[\Delta\cap \mathcal I^\infty(G)]$, i.e., $\gamma/G$ is solution of $\Delta_G\coloneqq \Pi^\infty[\Delta\cap \mathcal I^\infty(G)]$. \textit{Vice-versa}, if $\gamma/G$ is a solution of $\Delta_G$, then $\gamma$ is a solution of $\Delta$, as it can be seen inverting $\Pi$ (see Ref.~\cite{olver} about the details of this inversion procedure). The new manifold $\Delta_G$ corresponds in particular to the reduced equation we were looking for.

However, once the reduction has been performed, it will be useful if a choice of variables exists in such a way that we are able to express a solution of the original problem in the form $u=f(x)$. This is possible if we require an additional property for $G$, namely transversality. Transversality guarantees that the invariant solution under the action of $G$ are actually in the form $u=f(x)$ locally. If the algebra $\mathfrak g$ of $G$ is generated by
\begin{equation}
\mathbf v_k\coloneqq\sum_{i=1}^N\xi^i_k(x,u)\partial_{x_i}+\phi_k(x,u)\partial_u,\quad{k=1,\dots,r},
\end{equation}
where $r$ is the dimension of the algebra, the transversality condition in a neighbourhood of the point $(x,u)\in\M$ can be locally expressed as
\begin{equation}
\text{rank}\begin{pmatrix}
\xi^1_1(x,u)&\dots&\xi^N_1(x,u)\\
\vdots&&\vdots\\
\xi^1_r(x,u)&\dots&\xi^N_r(x,u)
\end{pmatrix}=\text{rank}\begin{pmatrix}\xi^1_1(x,u)&\dots&\xi^N_1(x,u)&\phi^1(x,u)\\
\vdots&&\vdots&\vdots\\
\xi^1_r(x,u)&\dots&\xi^N_r(x,u)&\phi^1(x,u)\end{pmatrix}=s,\label{transv}
\end{equation}
where $s$ is the dimension of the orbits of $G$. Due to the fact that the invariants are, by hypothesis, functionally independent, the Jacobian matrix
\begin{equation}
\begin{pmatrix}
\partial_{x_1}\zeta^1&\dots&\partial_{x_N}\zeta^1&\partial_u\zeta^1\\
\vdots&&\vdots&\vdots\\
\partial_{x_1}\zeta^{N-s}&\dots&\partial_{x_N}\zeta^{N-s}&\partial_u\zeta^{N-s}\\
\partial_{x_1}\upsilon &\dots&\partial_{x_N}\upsilon&\partial_u\upsilon
\end{pmatrix}
\end{equation}
has rank $N+1-s$ everywhere. Moreover, the transversality condition guarantees that the last column of the previous matrix cannot contain vanishing entries only \cite{olver}. This fact allow us to apply locally Dini's theorem and write $u$ and $N-s$ independent variables $\tilde x$ as a function of the $N-s+1$ invariants and of the $s$ remaining original independent variables.

\subsection{Example}\label{sec:3}
To illustrate our approach, let us now consider the following equation:
\begin{equation}
\mathcal E=\partial_t u-\sum_{i=1}^N\alpha_iu^{\beta_i}\rrfder[p_i]{0}{x_i}u=0,\label{genheat} \qquad\alpha_i,\beta_i\in\R\setminus\{0\},\quad p_i\in\R^+\setminus\N.
\end{equation}
The equation above has the structure of a fractional equation in which $t$ has the role of a time coordinate, whilst $x_i$, $i=1,\dots,N$, can be interpreted as space coordinates in an $N$--dimensional space. Observe that, in anomalous diffusion models, the time derivative is typically fractional, whilst the derivatives with respect to the other variables are of integer orders. Diffusion equations involving fractional derivatives with respect to position or velocity variables, however, are not uncommon. For example, equations of this type have been introduced in the study of diffusion of charged particle in a magnetic field \cite{Cechkin2002}, diffusion in heterogeneous disordered media \cite{Bologna2000,Lenzi2005}, or reaction--diffusion models \cite{Castillo2003}. Here we will discuss, however, Eq.~\eqref{genheat} mostly as a pedagogical example.

In Eq.~\eqref{genheat}, there are no mixed operators and only one integro-differential operator appears for each variable. We will solve the problem of determining symmetries for this equation under the hypothesis that fractional derivatives of different order are independent. As specified above, this working ansatz is fullfilled when the coefficients of the fractional derivatives in the determining equation depend at most on $x$, $u$ and a finite number of derivatives of $u$. Following the approach described above, and using as shorthand notation $x\equiv (x_1,\dots,x_N)$, we search for a Lie field of the form
\begin{equation}
\mathbf v=\tau(t,x,u)\frac{\partial}{\partial t}+\sum_{i=1}^N\xi^i(t,x,u)\frac{\partial}{\partial x_i}+\phi(t,x,u)\frac{\partial}{\partial u}.
\end{equation}
The determining equation is obtained applying Eq.~\eqref{genprolongation}
\begin{equation}
\calcin{\left\{\phi^{1}_{t}-\sum_{i=1}^N\left[\alpha_iu^{\beta_i}\phi^{p_i}_{x_i}+\alpha_i\beta_iu^{\beta_i-1}\phi\rrfder[p_i]{0}{x_i} u\right]\right\}}{\mathcal E=0}=0,
\end{equation}
In particular, here we have
\begin{multline}\label{phi0pgen}
\phi^p_{x_i}=\rrfder[p]{0}{x_i}\phi+\rrfder[p]{0}{x_i}u\left(\frac{\partial \phi}{\partial u}-p\frac{\dd\xi^i}{\dd x_i}\right)-u\rrfder[p]{0}{x_i}\frac{\partial \phi}{\partial u}-\sum_{n=1}^\infty\binom{p}{n}\frac{\dd^n\tau}{\dd x_i^n}\frac{\partial}{\partial t}\rrfder[p-n]{0}{x_i}u\\
+\sum\limits_{n=1}^\infty\left\{\left[\binom{p}{n}\frac{\partial^{n+1} \phi}{\partial x_i^n\partial u}-\binom{p}{n+1}\frac{\dd^{n+1}\xi^i}{\dd x_i^{n+1}}\right]\rrfder[p-n]{0}{x_i}u-\sum_{j\neq i}\binom{p}{n}\frac{\dd^n\xi^j}{\dd x_i^n}\frac{\partial}{\partial x_j}\rrfder[p-n]{0}{x_i}u\right\}\\+\sum\limits_{n=2}^\infty\sum\limits_{m=2}^n\sum\limits_{k=2}^m\sum\limits_{r=0}^{k-1}\binom{p}{n}\binom{n}{m}\binom{k}{r}\frac{{x_i}^{n-p}(-u)^r}{k!\Gamma(n+1-p)}\frac{\dd^m u^{k-r}}{\dd x_i^m}\frac{\partial^{n-m+k}\phi}{\partial x_i^{n-m}\partial u^k}.\end{multline}
After some manipulations, we get the following determining equation
\begin{multline}
\frac{\partial\phi}{\partial t}
=\sum_{i=1}^N\left\{\frac{\dd\xi^i}{\dd t}\frac{\partial u}{\partial x_i}+\left[\left(\frac{\dd\tau}{\dd t}-\frac{\partial\phi}{\partial u}\right)\alpha_iu^{\beta_i}
+\phi \alpha_i\beta_iu^{\beta_i-1}\right]\rrfder[p_i]{0}{x_i}u\right.\\
+\alpha_iu^{\beta_i}\sum_{n=0}^\infty\left[\binom{p_i}{n}\frac{\partial^{n+1}\phi}{\partial x_i^n\partial u}-\binom{p_i}{n+1}\frac{\dd^{n+1}\xi^i}{\dd x_i^{n+1}}\right]\rrfder[p_i-n]{0}{x_i}u\\
+\alpha_iu^{\beta_i}\left[\rrfder[p_i]{x_i}{0}\phi-u\rrfder[p_i]{0}{x_i}\frac{\partial\phi}{\partial u}-\sum_{n=1}^\infty\binom{p_i}{n}\left(\frac{\dd^n\tau}{\dd x_i^n}\frac{\partial}{\partial t}\rrfder[p_i-n]{0}{x_i}u-\sum_{j\neq i}\frac{\dd^n\xi^j}{\dd x_i^n}\frac{\partial}{\partial x_j}\rrfder[p_i-n]{0}{x_i}u\right)\right]\\
\left.+ \alpha_iu^{\beta_i}\sum\limits_{n=2}^\infty\sum\limits_{m=2}^n \sum\limits_{k=2}^m\sum\limits_{r=0}^{k-1}\binom{p}{n}\binom{n}{m}\binom{k}{r}\frac{{x_i}^{n-p}(-u)^r}{k!\Gamma(n+1-p)}\frac{\dd^m u^{k-r}}{\dd x_i^m}\frac{\partial^{n-m+k}\phi}{\partial x_i^{n-m}\partial u^k}\right\}.
\end{multline}
From the previous equation we deduce the relations
\begin{subequations}\label{45}\begin{align}
u^{\beta_i}\rrfder[p_i-n]{0}{x_i}u\to &\frac{n+1}{p_i-n}\frac{\partial^{n+1}\phi}{\partial x_i^n\partial u}-\frac{\dd^{n+1}\xi^i}{\dd x_i^{n+1}}=0,\quad \forall i,\quad \forall n\in\N\\
u^{\beta_i-1}\rrfder[p_i]{0}{x_i}u\to &\left(\frac{\dd\tau}{\dd t}-p_i\frac{\dd \xi^i}{\dd x_i}\right) u+\phi \beta_i=0,\quad \forall i,\\
\partial_t\rrfder[p_i-n]{0}{x_i}u\to&\frac{\dd^{n}\tau}{\dd x_i^{n}}=0\Rightarrow\tau\equiv \tau(t),\quad \forall i,\quad \forall n\in\N\\
\partial_{x_j}\rrfder[p_i-n]{0}{x_i}u\to&\frac{\dd^{n}\xi^j}{\dd x_i^{n}}=0\Rightarrow\xi^i\equiv \xi^i(t,x_i),\quad i\neq j,\quad \forall i,\quad \forall n\in\N\\
\partial_t u\to &\frac{\dd\xi^i}{\dd t}=0\Rightarrow\xi^i\equiv \xi^i(x_i),\quad\forall i,\quad \forall n\in\N.
\end{align}\end{subequations}
Imposing the additional constraint $\calcin{\xi^j}{x_j=0}=0$, $j=1,\dots,N$, the equations above have the following solutions
\begin{subequations}\begin{align}
\tau(t)&=c_1 t+c_0,\\
\xi^i(x_i)&=\frac{c_1-\beta_i c_2}{p_i} x_i,\quad i=1,\dots,N,\\
\phi&=c_2u,
\end{align}\end{subequations}
depending on the free real constants $c_0$, $c_1$ and $c_2$. We obtain therefore three symmetry generators,
\begin{subequations}
\begin{align}
\mathbf v_1&=\partial_{t},\\
\mathbf v_2&= t\partial_{t}+\sum_{i=1}^N\frac{x_i}{p_i}\partial_{x_i},\\
\mathbf v_3&= u\partial_u-\sum_{i=1}^N\frac{\beta_i x_i}{p_i}\partial_{x_i}.
\end{align}
\end{subequations}

The first generator expresses the trivial invariance of Eq.~\eqref{genheat} under translations in the variable $t$. The corresponding reduced equation can be obtained by looking for a stationary solution $u=\upsilon(x)$.

The generators $\mathbf v_2$ and $\mathbf v_3$ express the invariance of the equation under scaling transformations in the dependent and independent variables. To perform the reduction, in both cases, we will use the following simple identity, proved for example in Ref.~\cite{BL},
\begin{equation}\label{riscala}
\RRfder[p]{0}{t}u\left(\lambda^\alpha t\right)=\lambda^{\alpha p}\RRfder[p]{0}{\tau}u(\tau),\quad \tau\coloneqq\lambda^\alpha t,\quad\lambda,\alpha\in\R^+.
\end{equation}
Let us start analyze the generator $\mathbf v_2$. By means of the equation $\mathbf v_2 z(t,x,u)=0$,  we can determine $N$ invariants of the form
\begin{equation}
z_i= \frac{x_i}{t^\frac{1}{p_i}},\quad i=1,\dots,N,
\end{equation}
and search for a solution of Eq.~\eqref{genheat} of the type \begin{equation}u(t,x)\equiv \upsilon\left(\frac{x_1}{t^\frac{1}{p_i}},\dots,\frac{x_N}{t^\frac{1}{p_i}}\right).\end{equation} Using Eq.~\eqref{riscala}, the reduced equation for $\upsilon=\upsilon(z_1,\dots,z_N)$ is easily obtained in terms of the $N$ new variables $\{z_i\}_i$ as
\begin{equation}\label{reducedgenheat}
\sum_{i=1}^N\left[\frac{z_i}{p_i}\frac{\partial}{\partial z_i}+\alpha_i \upsilon^{\beta_i}\rrfder[p_i]{0}{z_i}\right]\upsilon=0.
\end{equation}

Similarly, the symmetry generator $\mathbf v_3$ provides the new set of variables
\begin{equation}
\upsilon=x_N^{-\frac{p_1}{\beta_1}}u,\qquad z_i=x_{i}x_N^{-\frac{\beta_{i} p_1}{\beta_1p_{i}}}\quad i=1,\dots,N-1.
\end{equation}
From them, we compute the reduced equation for $\upsilon=\upsilon(t,z_1,\dots,z_{N-1})$. We use again Eq.~\eqref{riscala}, obtaining
\begin{equation}
\sum_{i=1}^{N-1}\alpha_i u^{\beta_i}\rrfder[p_i]{0}{x_i}u=x^{\frac{p_1}{\beta_1}}\sum_{i=1}^{N-1}\alpha_i \upsilon^{\beta_i}\rrfder[p_i]{0}{z_i}\upsilon.
\end{equation}
Before proceeding further, we will introduce a new auxiliary function. Given a function of $n+1$ variables $f(t,x_1,\dots,x_n)\equiv f\left(x_1t^{-\mu_1},\dots,x_n t^{-\mu_n}\right)$, $\mu_i>0$ $\forall i$, and denoting by $\chi_i\coloneqq x_it^{-\mu_i}$, for a certain $\nu\in\R$ we can write
\begin{multline}
\rrfder[p]{0}{t}\left[t^\nu f(\{\chi_j\})\right]=\frac{\partial^{[p]+1}}{\partial t^{[p]+1}}\left[\int_0^t \frac{(t-\tau)^{[p]-p}\tau^\nu}{\Gamma([p]+1-p)}f\left(x_1\tau^{-\mu_1},\dots,x_n \tau^{-\mu_n}\right)\dd\tau\right]\\
=\frac{\partial^{[p]+1}}{\partial t^{[p]+1}}\left[t^{\nu+[p]-p+1}\mathsf K^{[p]-p,\nu}_{\left\{\mu_i\right\}_i}[f](\{x_j t^{-\mu_j}\})\right]\\
=t^{\nu-p}\prod_{k=0}^{[p]}\left(1+\nu-p+k-\sum_{i=1}^n\mu_i\chi_i\frac{\partial}{\partial \chi_i}\right)\mathsf K^{[p]-p,\nu}_{\left\{\mu_i\right\}_i}[f](\{\chi_j\})\eqqcolon t^{\nu-p}\mathsf P^{p,\nu}_{\left\{\mu_i\right\}_i}[f](\{\chi_j\}).
\end{multline}
The operator $\mathsf K^{p,\nu}_{\{\mu_i\}_i}[f]$ appearing in the relations above is a generalization of the Erd\'elyi--Kober fractional integral operator \cite{Erdelyi,Kober} and has the form
\begin{equation}
\mathsf K^{p,\nu}_{\{\mu_i\}_i}[f](\{z_i\})\coloneqq \int_0^1 \frac{(1-\tau)^{p}\tau^\nu}{\Gamma(p+1)}f\left(z_1\tau^{-\mu_1},\dots,z_n \tau^{-\mu_n}\right)\dd\tau.
\end{equation}
It follows that
\begin{equation}
\alpha_1u^{\beta_1}\rrfder[p_1]{0}{x_1}u=x^{\frac{p_1}{\beta_1}}\alpha_1\upsilon^{\beta_1}\mathsf P^{p,-\frac{p_1}{\beta_1}}_{\left\{\frac{\beta_ip_1}{\beta_1p_i}\right\}_i}[\upsilon](\{z_j\}).
\end{equation}
The reduced equation is therefore
\begin{equation}\label{reducedgenheat2}
\partial_t \upsilon=\alpha_1\upsilon^{\beta_1}\mathsf P^{p,-\frac{p_1}{\beta_1}}_{\left\{\frac{\beta_ip_1}{\beta_1p_i}\right\}_i}[\upsilon]+\sum_{i=1}^{N-1}\alpha_i \upsilon^{\beta_i} \rrfder[p_i]{0}{z_i}\upsilon.
\end{equation}

To conclude, we observe that the appearance of the Erd\'ely--Kober operator, or Erd\'elyi--Kober--type operators, after a symmetry reduction of a \textsc{fpde} is certainly not new in the literature \cite{BL}. Indeed, Erd\'elyi--Kober operators commonly appear in relaxation and oscillation models \cite{Concezzi2015}. In Ref.~\cite{Sahadevan2012}, a symmetry analysis of time-fractional Burgers and Korteweg--de Vries equations provided reduced equations involving this class of operators.

\section*{Acknowledgments}
G.\,S. acknowledges the financial support of the John Templeton Foundation. The research of P.\,T. has been partly supported by the research project FIS2015-63966, MINECO, Spain, and by the ICMAT Severo Ochoa project SEV-2015-0554 (MINECO).
\bibliographystyle{plain}
\bibliography{Biblio.bib}
\end{document}